\documentclass[11pt]{article}
\usepackage{amsmath}
\usepackage{latexsym}
\usepackage{exscale}
\def\cH{{\cal H}}
\def\cS{{\cal S}}
\newcommand{\set}[1]{\left\{#1\right\}}
\newcommand{\bin}[2]{\binom{#1}{#2}}
\newcommand{\rdown}[1]{\lfloor #1 \rfloor}
\newcommand{\brac}[1]{\left(#1\right)}
\newcommand{\bfrac}[2]{\brac{\frac{#1}{#2}}}
\newcommand{\whp}{{\bf whp}}
\newcommand{\Whp}{With high probability}

\def\cH{{\cal H}}

\def\cE{{\cal E}}
\def\cEk{{\cal E}^{(t)}}
\def\cA{{\cal A}}
\def\cB{{\cal B}}

\def\om{\omega}
\newtheorem{theorem}{Theorem}
\newtheorem{lemma}[theorem]{Lemma}

\newtheorem{proposition}[theorem]{Proposition}
\newtheorem{scholium}[theorem]{Claim}

\def\D{\Delta}

\newcommand{\proofstart}{{\bf Proof\hspace{2em}}}
\newcommand{\proofend}{\hfill\mbox{$\Box$}}

\newcommand{\nuall}{\nu_{{\rm all}}}
\newcommand{\nuemp}{\nu_\emptyset}
\newcommand{\nuf}{\nu_G}

\newcommand{\cM}{{\cal M}}
\def\t{\tau}

\def\three{t_3}
\def\tbigint{t_{\cap}}

\newcommand{\nuempA}{\nuemp^A}
\newcommand{\nuempB}{\nuemp^B}
\begin{document}
\makeatletter \title{On Randomly Generated Intersecting Hypergraphs II}
\author{Tom Bohman\thanks{ Supported in part by NSF grant DMS-0100400.}\\
\small Carnegie Mellon University,\\[-0.8ex]
\small \texttt{tbohman@andrew.cmu.edu}
\and Alan Frieze\thanks{Supported in part by NSF grant CCR-0200945.} \\
\small Carnegie Mellon University,\\[-0.8ex]
\small \texttt{alan@random.math.cmu.edu}
\and Ryan Martin\thanks{Contact author.  Supported in part by NSF VIGRE grant DMS-9819950 and NSA grant H98230-05-1-0257.} \\
\small Iowa State University,\\[-0.8ex]
\small \texttt{rymartin@iastate.edu} \and Mikl\'os
Ruszink\'o\thanks{Research was partially supported by OTKA Grants
T 046234,  T 038198 and by a NATO
Collaborative Linkage Grant}\\ 
\small Computer and Automation Research Institute\\ [-0.8ex]
\small Hungarian Academy of Sciences,\\ [-0.8ex] 
\small \texttt{ruszinko@sztaki.hu} \and Cliff
Smyth\thanks{Supported in part by NSF VIGRE grant DMS-9819950.  Part
of this research was done while the author was at Carnegie Mellon University.} \\ \small Massachusetts Institute of Technology,\\[-0.8ex]
\small \texttt{csmyth@math.mit.edu}}\date{} \vspace{-.1in} \maketitle

\begin{abstract}
Let $c$ be a positive constant. Suppose that  $r=o(n^{5/12})$ and
the members of $\binom{[n]}{r}$ are chosen sequentially at random
to form an intersecting hypergraph $\cH$. We show that
\whp\footnote{A sequence of events $\cE_1,\ldots,\cE_n,\ldots$ is
said to occur {\bf with high probability} (\whp) if
$\lim_{n\to\infty}\Pr(\cE_n)=1$.}\ $\cH$ consists of a simple
hypergraph $\cS$ of size $\Theta(r/n^{1/3})$, a distinguished
vertex $v$ and all $r$-sets which contain $v$ and meet every edge
of $\cS$. This is a continuation of the study of such random
intersecting systems started in~\cite{BCFRR} where the case
$r=O(n^{1/3})$ was considered. To obtain the stated result we
continue to investigate this question in the range $\omega
(n^{1/3})\le r \le o(n^{5/12})$.
\end{abstract}

\section{Introduction}
The study of random combinatorial structures has emerged as an
important component of Discrete Mathematics. The most intensely
studied area is that of random graphs \cite{B0}, \cite{JLR} and
many of the results of this area have been extended to hypergraphs
or set systems. There are (at least) two ways to study random
structures. One can set up a probability space and sample directly
from it or one can define a random process and study its outcome.
(One can argue that there is no formal difference between the two
approaches, but one cannot deny a qualitative difference between
them). As examples of the first approach, one can study the
distribution of the number of triangles in the random graph
$G_{n,p}$, or one can estimate the probability that the random
graph $G_{n,1/2}$ is triangle-free \cite{trianglefree}. As an
example of the second approach, one can consider a process where
we add random edges, but avoid introducing triangles \cite{ESW}.
In this paper, which is a continuation of \cite{BCFRR}, we follow
this idea and consider intersecting hypergraphs.

A {\bf hypergraph} is a family of subsets of a given ground set.
 The subsets are called {\bf edges}.  The {\bf degree} of a vertex
is the number of edges that contain it.  When we deal with more
than one hypergraph, the appropriate edges and degrees should be
clear from the context.  In the proofs, we will refer to auxiliary
{\bf graphs}, which are simple graphs in the ordinary sense.  That
is, $2$-uniform hypergraph.  An {\bf intersecting hypergraph} is
one in which each pair of edges has a non-empty intersection.
Here, we consider {\bf $r$-uniform hypergraphs} which are those
for which all edges contain $r$ vertices.

The motivating idea for this paper is the classical
Erd\H{o}s-Ko-Rado theorem \cite{EKR} which states that a maximum
size $r$-uniform intersecting hypergraph  has $\binom{n-1}{r-1}$
edges if $r\leq n/2$ and $\binom{n}{r}$ edges if $r>n/2$.
Furthermore, for $r<n/2$ any maximum-sized family must have the
property that all edges contain a common vertex.

In the last four decades this theorem has attracted the attention
of many researchers and it has been generalized in many ways. It
is worth mentioning, for example, the famous conjecture of Frankl
on the structure of maximum $t$-intersecting families in a certain
range of $n(t,r)$ which was investigated by Frankl and F\"uredi
\cite{FF} and completely solved only a few years ago by Ahlswede
and Khachatrian \cite{AK}. Another type of generalization can be
found in \cite{BFRT}.

The first attempt to `randomize' this topic was given by Fishburn,
Frankl, Freed, Lagarias and Odlyzko~\cite{FFF}.  Also note that
other random hypergraph structures were considered already by
R\'enyi e.g., in~\cite{R}, he identified the anti-chain threshold.
The paper~\cite{BCFRR} began the study of randomly generated
intersecting systems.  More precisely, the edges were taken
on-line; that is, one at a time, ensuring that at each stage, the
resulting hypergraph remained intersecting.  I.e., we considered
the following random process:

{\sc choose random intersecting system}\\
Choose $E_1\in \binom{[n]}{r}$.  Given
$\cH_t:=\{E_1,\ldots,E_t\}$, let $\cA(\cH_t)=\{E\in
\binom{[n]}{r}:\;E\notin \cH_t$ and $E\cap E_\tau\neq\emptyset$
for $1\leq \tau\leq t\}$. Choose $E_{t+1}$ uniformly at random
from $\cA(\cH_t)$. The procedure halts when $\cA(\cH_t)=\emptyset$
and $\cH=\cH_t$ is then output by the procedure.

It should be made clear that sets are chosen {\em without} replacement.

The paper~\cite{BCFRR} studied the case where $r=O(n^{1/3})$. The
main result of that paper was the following theorem which
determines the threshold for the event that edges chosen randomly
online to form an intersecting hypergraph will attain the
Erd\H{o}s-Ko-Rado bound:

We say that $\cal H$ {\bf fixes} $x$ if every member of
$\cal H$ contains $x$.

\begin{theorem}\label{main}\cite{BCFRR}
Let $\cE_{n,r}$ be the event that
$\left\{|\cH|=\binom{n-1}{r-1}\right\}$.  For $r<n/2$, this is
equivalent to $\cal F$ fixing some $x\in [n]$.  Then if
$r=c_nn^{1/3}<n/2$,
$$ \lim_{n\to\infty}\Pr(\cE_{n,r})=
   \begin{cases} 1 & c_n\to 0 \\
                 \frac{1}{1+c^3} & c_n\to c \\
                 0 & c_n\to\infty
   \end{cases}.$$
\end{theorem}
For $k\geq 3$ let
$$t_k=\min\{t:\;\D(\cH_t)=k\}$$
where for $v\in [n]$, $\deg_t(v)$ is the degree of $v$ in $\cH_t$ and $\D(\cH_t)=\max_{v\in V}\set{\deg_t(v)}$. 

As $r$ grows beyond $n^{1/3}$ the structure of $\cH$ grows more
complex. In this paper we are able to analyze a small part of the
range where $r/n^{1/3}\to\infty$: Hypergraph
$\cH_t$ is {\bf simple} if $|E_i\cap
E_j|\leq 1$ for all $1\leq i\neq j\leq t$. In the range investigated, we show
that \whp\ $\cH_t$ remains simple and has maximum degree two for the first $\Theta(r/n^{1/3})$
edges. Then at time $t_3$ a vertex $v$ of degree
three is created, then $o(t_3)$ edges are added that avoid $v$
until, at time $t_4$, $v$ achieves degree four.  Then, \whp, the process
finishes with all the subsequent edges containing $v$ and intersecting
the $t_4-4$ edges of $\cH_{t_4}$ not containing $v$.

\begin{theorem}\label{th2}
   Suppose that $\om\to\infty$ and
   \begin{equation}\label{eqr}
      \om n^{1/3}\leq r\leq n^{5/12}/\om.
   \end{equation}
Then
\begin{enumerate}
\item $nt_4^3/(6r^3)$  converges to the exponential distribution with
  mean 1 i.e.
$$ \lim_{n\to\infty}\Pr(t_4\geq cr/n^{1/3})= e^{-c^3/6}.$$
\item $\cH_{t_4}$ is simple \whp.
\item $\cH_{t_4}$ has a unique vertex $v$ of degree 4 and no vertex of
  degree 3, \whp.
\item $\cH$ is the unique hypergraph consisting
   of all edges that (a) contain $v$ and (b) meet every edge of
   $\cH_{t_4}$ which does not contain $v$.
\end{enumerate}
    Furthermore, given $\cH_{t_4}$ satisfying the above,
   \begin{eqnarray}
      |\cH|&=&\sum_{i=0}^{t_4-4}(-1)^i\binom{t_4-4}{i}%
            \binom{n-1-ri+i(i-1)/2}{r-1}.\label{exactH}\\
&=&(1+o(1))\left(\frac{r^2}{n}\right)^{t_4-4}%
                    \bin{n-1}{r-1} \label{approxH}
   \end{eqnarray}
\end{theorem}

We break the proof of  Theorem~\ref{th2} into stages. In
Section~\ref{Ph1} we study the growth of $\cH$ until the time
$t_3$ when the first vertex of degree at
least three is created. As part of the analysis we show that \whp\
$\cH$ remains simple up to this point. Subsequent sections inch
along until the maximum degree is
$\D_0=\rdown{(r/n^{1/3})^{1/2}}$. At this time, \whp, $\cH$ is
still simple and $v$ is the unique vertex of degree more than two.
Section~\ref{Ph4} then shows that things finish as described in
the theorem. Section \ref{gridfun} provides estimates for some quantities
used in the proofs.

In this and the previous work \cite{BCFRR} we have tried to
examine the structure of a typical intersecting family, in particular
the one grown by a natural sequential process. Our proof requires
long and careful computation, but more than this, we needed strong
insight into the final goal to guide us. We hope to continue this
study for other ranges of $r$, but we do not at the moment see what
the typical structure is like for larger $r$.
\section{Proof of Theorem~\ref{th2}}

\subsection{The first degree-three vertex}
\label{Ph1}

In this section we study the growth of $\cH$ until the first time
$\three$ that $\cH$ contains a vertex of degree 3. We define
$$\tbigint=\min\{t:\;\exists s<t:\;|E_s\cap E_t|\geq 2\}$$
to be the first time that $\cH_t$ is not simple.

Our aim now is to show that the following behavior is typical: For
any constant $c>0$, the probability that $\three\geq cr/n^{1/3}$
falls off like $e^{-c^3/6}$ and that \whp\ at time $\three$,
$\cH_{\three}$ is simple.

So we split the possibilities into various events. We consider:
\begin{itemize}
   \item $\cA_t\stackrel{\rm def}{=}
   \{\D(\cH_t)\leq 2\}$
   \item $\cA_t^v\stackrel{\rm def}{=}
   \left\{\deg_t(w)\leq 2,
   \forall w\neq v\mbox{ and }\deg_t(v)\geq 3\right\}$
   \item $\cB_t\stackrel{\rm def}{=}
   \left\{|E_i\cap E_j|=1, \forall 1\leq i<j\leq t\right\}$
\end{itemize}

In words, $\cA_t$ is the event that the hypergraph $\cH_t$ has
maximum degree 2, $\cA_t^v$ is the event that $\deg_t(v)\geq 3$
and that the maximum degree of the remaining vertices $w\neq v$ is
2 and $\cB_t$ is the event that $\cH_t$ is simple.

We partition the possibilities for $\cH_t$ into three disjoint
events and later we will express them in terms of the
$\cA_i,\cA^v_i,\cB_i$.
\begin{eqnarray*}
   \cEk_1 & = & \{\three>t\}\wedge\{\tbigint>t\}\\
   \cEk_2 & = & \{\three\leq t\}\wedge\{\three<\tbigint\}%
                \wedge\{\exists v : \cA_{\three}^v\}\\
   \cE_3 & = & \bigcup_t\brac{\overline{\cEk_1}\cap \overline{\cEk_2}}.
\end{eqnarray*}

\begin{lemma}\label{t3main}
   Assume that
\begin{equation}\label{uppert}
t\leq \frac{Kr}{n^{1/3}}
\end{equation}
where $K$ is some arbitrarily large positive constant.

   \begin{description}
      \item[(a)] $\Pr(\cEk_1)=%
      \exp\left\{-\frac{t^3n}{6r^3}\right\}+O(t^5n^2/r^6+t^2r^2/n+t^2n/r^3)$.
      \item[(b)] $\Pr(\cEk_2)=%
      1-\exp\left\{-\frac{t^3n}{6r^3}\right\}+O(t^5n^2/r^6+t^2r^2/n+t^2n/r^3)$.
      \item[(c)] $\Pr\brac{\cE_3}=o(1)$.
   \end{description}
Here the hidden constant depends on $K$.
\end{lemma}
\proofstart We write
\begin{eqnarray}
   \Pr(\cEk_1) & = & \prod_{i=1}^{t-1}%
   \Pr\left(\cA_{i+1}\wedge\cB_{i+1}\mid
   \cA_i\wedge\cB_i\right) \label{sum1} \\
   \Pr(\cEk_2) & = & \sum_{i=1}^{t-1}%
   \Pr\left(\left\{\bigvee_v \cA_{i+1}^v\right\}\wedge%
   \cB_{i+1}\mid \cA_i\wedge\cB_i\right)\Pr(\cA_i\wedge\cB_i)%
   \label{sum2}
\end{eqnarray}
Note $\Pr(\cA_1\wedge\cB_1)=1$.

Consider the following quantities:
\begin{itemize}
   \item $\nuall(t)$ is the number of $r$-sets that
   intersect every edge of $\cH_t$.
   \item Assume that $\cH_t$ is a simple hypergraph and that
   $E_i\cap E_j=\{x_{i,j}\}$ for $1\leq i<j\leq t$.  Let
   $G=(S,F)$ be a graph with $s$ vertices and $f$ edges where
   $S\subseteq [t]$.
   \begin{itemize}
      \item $\nuf^{*}(t)$ is the number of $r$-sets $E$ such
      that $x_{i,j}\in E$ iff $(i,j)\in F$.
      \item $\nuf(t)$ is the number of $r$-sets $E$ such that
      $x_{i,j}\in E$ iff $(i,j)\in F$ and which meet each
      $E_i,\,i\notin S$ in exactly one vertex.
      \item $\cM$ denotes the set of $G$ which are matchings.
   \end{itemize}
   \item As we will see, the dominant term $\nuemp(t)$ is the number of $r$-sets that (i) intersect every
   edge of $\cH_t$, (ii) keep $\cH_{t+1}$ simple, given that
   $\cH_t$ is simple and (iii) keep $\D(\cH_{t+1})\leq 2$, given
   $\D(\cH_t)\leq 2$. \\ (This is the case $f=0$.)
\end{itemize}

Note that $s\leq t$ and $f\leq
\binom{s}{2}$ and $f=0$ implies $s=0$.
\begin{eqnarray}
   \nuf(t) & = & (r-t+1)^{t-s}\binom{n-t(r-t+1)-\binom{t}{2}}%
                                    {r-t+s-f} \nonumber\\
           & = & (1+O(tr^2/n))r^{t-s}\frac{n^{r-t+s-f}}%
                                          {(r-t+s-f)!} \label{empty}\\
           & = & (1+O(tr^2/n))\bfrac{n}{r^2}^s%
                              \bfrac{r}{n}^f\nuemp(t).\nonumber
\end{eqnarray}

We have taken care to remove extraneous error terms in our ``big O'' notation, based on the bound \eqref{uppert} that we have given for $t$.\\
Furthermore, we obtain an expression
$$\nuemp(t)=(1+O(tr^2/n))n^{r-t}r^{t}/(r-t)!$$
by putting $s=f=0$ into
\eqref{empty}.

Continuing our estimates,
\begin{equation}\label{qwer}
   \nuf^{*}(t)-\nuf(t)\leq
   \sum_{p=t-s+1}^{r-f}\ \binom{n}{r-f-p}%
   \sum_{\substack{a_j\geq 1,j\in [t-s] \\ \sum a_j=p}}
   \ \prod\limits_{j=1}^{t-s}\binom{r}{a_j}.
\end{equation}
{\bf Explanation:} The integer $p$ denotes the number of elements
of $E$ that belong to $\cH_t$ but do not lie in an edge
corresponding to a vertex of $G$. The quantity $a_j$ is the number
of elements of $E$ which lie in the $j^{\rm th}$ such edge. Having
chosen these $p$ elements and the $f$ elements of $\cH_t$
corresponding to the edges of $G$, we have at most
$\binom{n}{r-f-p}$ choices for the remaining elements of $E$.

Increasing $p$ by 1 reduces the first binomial coefficient of
\eqref{qwer} by a factor of $\sim r/n$. We can get all choices
satisfying $a_1+\cdots+a_{t}=p+1$ by choosing
$a_1'+\cdots+a_{t}'=p$, choosing a $j$ and adding one to $a_j'$.
Thus the number of choices increases by at most $tr$ and so the
ratio of $p+1$ terms to $p$ terms is $O(tr^2/n)$.

So,
\begin{eqnarray*}
   \nuf^{*}(t)-\nuf(t)
   & \leq & (1+O(tr^2/n))\frac{tr^{t-s+1}}{2}%
                         \binom{n}{r-f-t+s-1} \\
   & \leq & (1+O(tr^2/n))\frac{tr^{t-s+1}}{2}%
                         \frac{n^{r-t+s-f-1}}{(r-t+s-f-1)!} \\
   & =    & (1+O(tr^2/n))\frac{tr^2}{2n}\nuf(t) \\
   & \leq & (1+O(tr^2/n))\frac{tr^2}{2n}\bfrac{n}{r^2}^{s}%
                         \bfrac{r}{n}^{f}\nuemp(t).
\end{eqnarray*}
Applying Proposition \ref{Gfun} (see Section \ref{appA}) with
$x=\frac{r}{n}$, $y=\frac{n}{r^2}$  we get:
\begin{eqnarray*}
   \nuemp(t)^{-1}\sum_{\substack{G\in\cM\\f=1}}\nu_{G}(t)
   & = & (1+O(tr^2/n))\binom{t}{2}\frac{n}{r^3} \\
   \nuemp(t)^{-1}\sum_{\substack{G\in\cM\\f\geq 2}}\nu_{G}(t)
   & = & O(t^4n^2/r^6) \\
   \nuemp(t)^{-1}\left(\nuemp^*(t)-\nuemp(t)\right)
   & = & O(tr^2/n) \\
   \nuemp(t)^{-1}\sum_{\substack{G\in\cM\\f\geq 1}}%
                 (\nu^*_{G}(t)-\nu_{G}(t))
   & = & O(t^3/r) \\
   \frac{\nu_{rest}(t)}{\nuemp(t)} & = & O(t^4n/r^4).
\end{eqnarray*}
The quantity
$$\nu_{rest}(t)=\sum_{\substack{G\notin \cM\\f\neq 0}}\nu^*_G(t)$$
 accounts for every possibility not specifically
mentioned in the previous four ratios.

We remark for future use that this implies that if \eqref{uppert} holds then
\begin{equation}
   \nuall(t)=\left(1+O(tr^2/n+t^2n/r^3)\right)\nuemp(t)%
   =\left(1+O(tr^2/n+t^2n/r^3)\right)r^t\binom{n}{r-t} .
   \label{future}
\end{equation}

Continuing, it further follows that
\begin{eqnarray}
   \Pr\left(\cA_{t+1}\wedge\cB_{t+1}\mid\cA_t\wedge\cB_t\right)
   & = & \frac{\nuemp(t)}{\nuall(t)} \nonumber \\
   & = & 1-\frac{t^2n}{2r^3}+O(t^4n^2/r^6+tr^2/n+tn/r^3) . \\
   \Pr\left(\left\{\bigvee_v\cA_{t+1}^v\right\}\wedge
   \cB_{t+1}\mid\cA_t\wedge\cB_t\right)
   & = & \frac{\nuemp(t)}{\nuall(t)}%
   \sum_{\substack{G\in\cM\\f=1}}\frac{\nuf(t)}{\nuemp(t)}
   \nonumber \\
   & = & \frac{t^2n}{2r^3}+O(t^4n^2/r^6+tr^2/n+tn/r^3) . \\
   \Pr\left(\cA_{t+1}\wedge\overline{\cB_{t+1}}\mid
   \cA_t\wedge\cB_t\right)
   & = & \frac{\nuemp(t)}{\nuall(t)}%
   \frac{\nuemp^{*}(t)-\nuemp(t)}{\nuemp(t)} \nonumber \\
   & = & O(tr^2/n) \label{sum3a} \\
   & & \nonumber \\
   \Pr\left(\bigwedge_v\overline{\cA_{t+1}^v}\wedge\cB_{t+1}\mid
   \cA_t\wedge\cB_t\right) & = &
   \frac{\nuemp(t)}{\nuall(t)}
   \sum_{\substack{G\in\cM \\ f\geq 2}}\frac{\nuf(t)}{\nuemp(t)}
   \nonumber \\
   & = & O(t^4n^2/r^6) \label{sum4a} \\
   & & \nonumber \\
   \Pr\left(\overline{\cA_{t+1}}\wedge\overline{\cB_{t+1}}\mid
   \cA_t\wedge\cB_t\right)
   & = & \frac{\nuemp(t)}{\nuall(t)}
         \left[\frac{\nu_{\rm rest}^{*}(t)}{\nuemp(t)}%
               +\sum_{\substack{ G\in\cM \\f\geq 1}}%
                \frac{\nuf^{*}(t)-\nuf(t)}{\nuemp(t)}\right]
   \nonumber \\
   & = & O(t^3/r). \label{sum5a}
\end{eqnarray}

Thus, if \eqref{uppert} holds,
\begin{eqnarray*}
   \Pr(\cEk_1)= \Pr(\cA_t\wedge\cB_t)
   & = & \prod_{j=1}^{t-1}\brac{1-\frac{j^2n}{2r^3}%
         +O(t^4n^2/r^6+tr^2/n+tn/r^3)}\\
   & = & \prod_{j=1}^{t-1}\exp\left\{-\frac{j^2n}{2r^3}+O(t^4n^2/r^6+tr^2/n+tn/r^3)\right\}\\
   & = & \exp\left\{-\frac{t^3n}{6r^3}\right\}%
         +O(t^5n^2/r^6+t^2r^2/n+t^2n/r^3),
\end{eqnarray*}
proving (a).

For $\cEk_2$ we have
\begin{eqnarray*}
   \Pr(\cEk_2)
   & = & \sum_{j=1}^t\frac{j^2n}{2r^3}e^{-j^3n/6r^2}%
         +O(t^5n^2/r^6+t^2r^2/n+t^2n/r^3) \\
   & = & 1-\exp\left\{-\frac{t^3n}{6r^3}\right\}%
         +O(t^5n^2/r^6+t^2r^2/n+t^2n/r^3)
\end{eqnarray*}
and (b) follows.

To prove (c) expand
$$\overline{\cEk_1}\cap \overline{\cEk_2}=(\{t_3\leq t\}\vee\{\tbigint\leq t\})\cap(\{t_3>t\}\vee\{t_3\geq \tbigint\}\vee\{\not\exists v:\cA_{\three}^v\}).$$
We see then that
$$\cE_3=\{t_3\geq \tbigint\}\vee\{\not\exists v:\cA_{\three}^v\}).$$
Now let $t_K=Kr/n^{1/3}$. From parts (a),(b) we see that
\begin{eqnarray*}
\Pr(t_3\geq \tbigint)&\leq& \Pr(t_3\geq t_K)+\Pr(\tbigint\leq t_K)\\
&\leq&\Pr(\overline{\cE^{(t_K-1)}_2})+\Pr(\overline{\cE^{(t_K)}_1})\\
&=&o(1).\\
\Pr(\not\exists v:\cA_{\three}^v)&\leq&\Pr(\overline{\cE^{(t_K)}_2})\\
&=&e^{-K^3/6}+o(1).
\end{eqnarray*}
We can make $K$ as large as we like and (c) follows.
\proofend

The following summarizes what we have proved so far:
\begin{scholium}\label{schol1}
   \Whp\ there $v\in V(\cH_{t_3})$ such that
   $\cH_{t_3}$ is simple and $v$ is the unique vertex of degree
   more than 2. Furthermore,
   $$ \lim_{n\to\infty}\Pr\brac{t_3\geq \frac{\alpha r}{n^{1/3}}}=e^{-\alpha^3/6}. $$
\end{scholium}
\proofstart
This follows from $\Pr(\cE^{(t_\alpha)}_1)=\Pr(t_3\geq t_\alpha)-o(1)$ where $t_\alpha=\frac{\alpha r}{n^{1/3}}$.
\proofend

\subsection{A Useful Lemma}\label{Ph2and3}

\newcommand{\tdel}{\t_\Delta}
In this section, we develop a lemma that can be used in each of
the next two sections.  It analyzes the behavior after $t_3$
when the maximum degree is still relatively small.  We will be
given $\cH_t$, an intersecting hypergraph that satisfies
$\cA_t^v\wedge\cB_t$ with maximum degree $\D\geq 3$.

Denote the following:
\begin{itemize}
   \item $\nuall^A(t)$ is the number of remaining $r$-sets
   that have a non-empty intersection with all edges of $\cH_t$ and
   contain $v$.

   \item $\nuall^B(t)$ is the number of remaining $r$-sets
   that have a non-empty intersection with all edges of $\cH_t$ but
   do not contain $v$.

   \item $\nuemp^A(t)$ is the number of $r$-sets that
   intersect the edges of $\cH_t$ containing $v$ in $v$ {\em only}
   and intersect the remaining $\tdel=t-\D$ edges in
   exactly one vertex, never creating any vertex of degree
   greater than $2$ (other than $v$).

   \item $\nuemp^B(t)$ is the number of $r$-sets that do not
   contain $v$, but do intersect each of the $t$ edges in
   exactly one vertex, also never creating any vertex of degree
   greater than $2$.
\end{itemize}

\begin{lemma} \label{smdel}
   Let $\cH_t$ be an intersecting hypergraph that satisfies
   $\cA_t^v\wedge\cB_t$ with maximum degree
   $\Delta$, $3\leq\D<\D_0$, where $\D_0=\lfloor
   (r/n^{1/3})^{1/2}\rfloor$.  Furthermore, assume that
\begin{equation}\label{Theta}
t=\Theta\bfrac{r}{n^{1/3}}.
\end{equation}
With the above
   notation,
\begin{eqnarray}
 \frac{\nuall^A(t)}{\nuemp^A(t)}&=&1+O\left(\frac{t^2n}{r^3}\right)\label{lem6a}\\
 \frac{\nuall^B(t)}{\nuemp^B(t)}&=&1+O\left(\frac{t^2n}{r^3}\right)\label{lem6b}\\
 \frac{\nuemp^B(t)}{\nuemp^A(t)}&=&\left(1+O\left(\frac{tr^2}{n}\right)\right)\frac{n}{r}\left(\frac{r^2}{n}\right)^\D\label{lem6c}
\end{eqnarray}
\end{lemma}

\proofstart

First we compute the two main expressions.
\begin{eqnarray}
   \nuemp^A(t) & = & (r-t+1)^{\tdel}%
   \binom{n-(rt-\binom{\tdel}{2}-\D \tdel-(\D-1))}{r-\tdel-1}
   \nonumber \\
   & = & \left(1+O(tr^2/n)\right)%
   \frac{r^{\tdel}n^{r-\tdel-1}}{(r-\tdel-1)!} \label{t4} \\
   \nuemp^B(t) & = & (r-t+1)^{\tdel}(r-\tdel-1)^\D%
   \binom{n-(rt-\binom{\tdel}{2}-\D \tdel-(\D-1))}{r-t} \nonumber \\
   & = & \left(1+O(tr^2/n)\right)\frac{r^tn^{r-t}}{(r-t)!}
   \nonumber \\
   & = & \left(1+O(tr^2/n)\right)\frac{n}{r}%
   \bfrac{r^2}{n}^\Delta\nuemp^A(t) \label{B}
\end{eqnarray}
Equation \eqref{lem6c} follows immediately.

Re-number the edges so that $v\in E_{\tdel+i},1\leq i\leq \D$.
Suppose that $E_i\cap E_j=\{x_{i,j}\}$ for $1\leq i<j\leq\tdel$
and that $E_i\cap E_{\tdel+j}=\{y_{i,j}\}$ for $1\leq i\leq\tdel$
and $1\leq j\leq \D$.

Let $G=(S,F)$ be a graph with $s\geq 2$ vertices and $f\geq 1$
edges where $S\subseteq [\tdel]$ and $F$ spans $S$. We define
$\nu_G^{A_*}(t)$ to be the number of $r$-sets $E$ such that (i)
$v\in E$, (ii) $x_{i,j}\in E$ iff $(i,j)\in F$. For this count we are dropping the
condition that the new edge intersects the old edges in exactly one
vertex.
\begin{eqnarray}
   \nu_G^{A_*}(t) & \leq &
   \sum_{p\geq\tdel-s}\binom{n}{r-1-f-p}%
   \sum_{\substack{a_j\geq 1, j\in[\t_\D-s] \\ \sum a_j=p}}%
   \prod_{j=1}^{\tdel-s}\binom{r}{a_j} \label{poiu} \\
   & = & \left(1+O(tr^2/n)\right)r^{\tdel-s}%
   \binom{n}{r-1-f-\tdel+s}, \nonumber
\end{eqnarray}
since increasing $p$ by 1 reduces the first binomial coefficient
in \eqref{poiu} by a factor of $\sim r/n$ and then we gain at most
a further $tr$ factor in the number of choices for the $a_j$.

So,
\begin{eqnarray*}
   \nu_G^{A_*}(t)  & \leq & \left(1+O(tr^2/n)\right)%
   r^{\tdel-s}\frac{n^{r-1-f-\tdel+s}}{(r-1-f-\tdel+s)!} \\
   & = & \left(1+O(tr^2/n+f^2/r)\right)\left(\frac{r}{n}\right)^f%
   \left(\frac{n}{r^2}\right)^s\nuemp^A(t).
\end{eqnarray*}

Finally, we compute
\begin{eqnarray*}
   \nuemp^{A_*}(t)-\nuemp^A(t) & \leq & \sum_{p\geq \tdel+1}\bin{n}{r-1-p}
   \sum_{\substack{a_j\geq 1, j\in[\t_\D-s] \\ \sum a_j=p}}%
   \prod_{j=1}^{\tdel}\binom{r}{a_j} \\
   & = & \left(1+O(tr^2/n)\right)\bin{n}{r-\tdel-2}\frac{\tdel}{2}
   r^{\tdel}(r-1) \\
   & = & \left(1+O(tr^2/n)\right)\frac{\tdel r^2}{2n}%
   \nuemp^A(t)
\end{eqnarray*}

So, from Proposition \ref{Gfun} of Section \ref{appA},
\begin{eqnarray*}
   \frac{\nuall^A(t)}{\nuempA(t)}-1 & \leq &
   \left(1+O(tr^2/n)\right)\left[%
   \sum_{\emptyset\neq G}\left(\frac{r}{n}\right)^f%
   \left(\frac{n}{r^2}\right)^s+\frac{\tdel r^2}{2n}\right]
   \nonumber \\
   & = & O(t^2n/r^3)
\end{eqnarray*}
proving \eqref{lem6a}.

Now for graph $G=(S,F)$ and bipartite graph $H\subseteq
([\tdel]\setminus S)\times [\Delta]$ define $\nu^{B_*}_{G,H}(t)$ to be
the number of $r$-sets $E$ such that (i) $v\notin E$, (ii)
$x_{i,j}\in E$ iff $(i,j)\in F$ and (iii) $y_{i,j}\in E$ iff
$(i,j)\in H$. Also let $u=u_H=|\{x:\;\exists y\ s.t.\ (x,y)\in
H\}|$ and $\ell=\ell_H=|\{y:\;\exists x\ s.t.\ (x,y)\in H\}|$ and $h=|E(H)|$.
Suppose either $G$ or $H$ is non-empty.
\begin{eqnarray*}
   \lefteqn{\nu^{B_*}_{G,H}(t)} \nonumber \\
   & \leq & \sum_{p\geq\tdel-s-u}\sum_{q\geq \Delta-\ell}%
   \binom{n}{r-f-h-p-q}%
   \sum_{\substack{a_j\geq 1, j\in[\tdel-s-u] \\ \sum a_j=p}}%
   \prod_{j=1}^{\tdel-s-u}\binom{r}{a_j}%
   \sum_{\substack{b_j\geq 1, j\in[\D-\ell] \\ \sum b_j=q}}%
   \prod_{j=1}^{\Delta-\ell}\binom{r}{b_j} \label{popo} \\
   & = & \left(1+O(tr^2/n)\right)\binom{n}{r-f-h-t+s+u+\ell}%
   \sum_{\substack{i_j\geq 1, j\in[\tdel-s-u] \\ \sum a_j=\tdel-s-u}}%
   \prod_{j=1}^{\tdel-s-u}\binom{r}{a_j}%
   \sum_{\substack{b_j\geq 1,j\in[\D-\ell] \\ \sum b_j=\Delta-\ell}}%
   \prod_{j=1}^{\Delta-\ell}\binom{r}{b_j} \nonumber \\
   & = & \left(1+O(tr^2/n)\right)r^{t-s-u-\ell}%
   \frac{n^{r-f-h-\tdel+s+u-\Delta+\ell}}{({r-f-h-t+s+u+\ell})!}
   \nonumber \\
   & = & \left(1+O(tr^2/n)\right)%
   \bfrac{r}{n}^{f+h}\bfrac{n}{r^2}^{s+u+\ell}\nuemp^B(t)
   \nonumber
\end{eqnarray*}

Finally,
\begin{eqnarray*}
   \nu^{B_*}_{\emptyset,\emptyset}(t)-\nuemp^B(t)
   & \leq & \sum_{p\geq\tdel+1}\sum_{q\geq \Delta+1}%
   \binom{n}{r-p-q}%
   \sum_{\substack{a_j\geq 1, j\in[\t_\D] \\ \sum a_j=p}}%
   \prod_{j=1}^{\tdel}\binom{r}{a_j}%
   \sum_{\substack{b_j\geq 1, j \in[\t_\D]\\ \sum b_j=q}}%
   \prod_{j=1}^{\Delta}\binom{r}{b_j} \\
   &   & +\sum_{q\geq \Delta+1}\binom{n}{r-\tdel-q}%
   \left(\frac{\tdel r^{\tdel}(r-1)}{2}\right)%
   \sum_{\substack{b_j\geq 1, j \in[\t_\D]\\ \sum b_j=q}}%
   \prod_{j=1}^{\Delta}\binom{r}{b_j} \\
   &   & +\sum_{p\geq\tdel+1}\binom{n}{r-p-\D}%
   \sum_{\substack{a_j\geq 1,  j \in[\t_\D] \\ \sum a_j=p}}%
   \prod_{j=1}^{\tdel}\binom{r}{a_j}%
   \left(\frac{\D r^\D(r-1)}{2}\right) \\
   & = & \left(1+O(tr^2/n)\right)\frac{tr^t(r-1)n^{r-t-1}}%
   {2(r-t-1)!} \\
   & = & \left(1+O(tr^2/n)\right)\frac{tr^2}{2n}\;\nuemp^B(t) .
\end{eqnarray*}

So,
$$ \frac{\nuall^B(t)}{\nuemp^B(t)}-1\leq
   \left(1+O(tr^2/n)\right)
   \left[\frac{tr^2}{2n}+\sum_G^*\sum_H^*%
         \bfrac{r}{n}^f\bfrac{n}{r^2}^s\bfrac{r}{n}^h%
         \bfrac{n}{r^2}^{u+\ell}\right] . $$
(The notation $\displaystyle \sum^*\sum^*$ means that the
summation avoids $(G,H)=(\emptyset,\emptyset)$)

Now from Propositions~\ref{Gfun} and~\ref{16} of
Section~\ref{appA}, ($\t=\tdel-s$, $\t\D
xy^2=\t\D(r/n)(n^2/r^4))\leq \tdel\D n/r^3=o(1)$),
\begin{align*}
   &\frac{\nuall^B(t)}{\nuemp^B(t)}-1\\
  &\leq\left(1+O(tr^2/n)\right)\left[\frac{tr^2}{2n}
       +\sum_{G\neq\emptyset}\left(\frac{r}{n}\right)^f
        \left(\frac{n}{r^2}\right)^s
        \left(1+(1+o(1))\frac{t^2n}{r^3}\right)
       +(1+o(1))\frac{\tdel\D n}{r^3}\right] \\
   &=O(t^2n/r^3)
\end{align*}
proving \eqref{lem6b}.
\proofend

\subsection{From $t_3$ to $t_4$}\label{Ph2}
Assume that $t_3=\Theta(r/n^{1/3})$. Let $t_4$ be the first time there is a vertex of degree 4.

It follows from Lemma~\ref{smdel} that as long as \eqref{Theta} holds and $\D\leq \D_0$,
\begin{eqnarray}
   \frac{\nuempA(t)+\nuempB(t)}{\nuall(t)}
   & = & 1-O(t^2n/r^3). \label{D} \\
   \frac{\nuempB(t)}{\nuempA(t)}
   & = & \left(1+O(tr^2/n)\right)\frac{r^5}{n^2},
   \qquad \D=3 . \label{E} \\
   \frac{\nuempB(t)}{\nuempA(t)}
   & = & O(r^7/n^3) . \qquad \D\geq 4. \label{F}
\end{eqnarray}

Let $t=t_3+\xi$ where $\xi\geq 0$ and $\xi=O(r^5/n^2)=o(r/n^{1/3})$. Note that
$r=o(n^{5/12})$ implies that $r^5/n^2=o(r/n^{1/3})$. Then
\eqref{D} implies
$$ \Pr(\cA_t^v\wedge\cB_t)=1-O(t^2r^2/n)=1-o(1) . $$
It follows from \eqref{D} and \eqref{E} that
\begin{eqnarray}
\Pr(t_4>t)&=&\prod_{\t=t_3+1}^t \frac{1-O(\t^2n/r^3)}{1+(1-O(\t r^2/n))n^2/r^5}\nonumber\\
&=&e^{O(\xi t_3^2n/r^3)}e^{O(t_3n/r^5)}\bfrac{r^5}{n^2+r^5}^{\xi}\nonumber\\
&=&e^{o(1)}\bfrac{r^5}{n^2+r^5}^{\xi}.\label{G}
\end{eqnarray}
If $r=o(n^{2/5})$ then the RHS of~\eqref{G} is $o(1)$ for all
$\xi\geq 1$ and we deduce that in this case $t_4=t_3+1$ \whp.

If $r=an^{2/5}$ where $a=a(n)\to\infty$ is allowed and $\xi=cr^5/n^2$
where $c>0$ is constant then
$$ \Pr(t_4> t)=e^{o(1)}(1+a^{-5})^{-ca^5}. $$
We deduce that
$$ \Pr(\neg(\cA_{t_4}^v\wedge\cB_{t_4}))\leq%
   o(1)+e^{o(1)}(1+a^{-5})^{-ca^5}.$$
Since $c$ can be made arbitrarily large, we deduce that
at time $t_4$ there is \whp, a unique vertex $v$ of degree $>2$.
The following summarizes what we have proved in this section:
\begin{scholium}\label{schol2}
   \Whp\ $t_4=(1+o(1))t_3$ and there exists $v\in V(\cH_{t_4})$
   such that $\cH_{t_4}$ is simple and $v$ is the unique vertex
   of degree more than two.
\end{scholium}
\proofend

\subsection{Finishing the proof}\label{Ph4}
Assume now that $\cA_{t_4}^v\wedge\cB_{t_4}$ holds and that
$t_4=(1+o(1))t_3$. The probability that in the next $\D_0-4$ steps
we either (i) add an edge not containing $v$ or (ii) that we make
a new vertex of degree 3 is at most
$O(\D_0(t_3^2n/r^3+r^7/n^3))=o(1)$.

Assume then that $\cA_{t_{\D_0}}^v\wedge\cB_{t_{\D_0}}$ holds.
Recall that $\D_0=\rdown{(r/n^{1/3})^{1/2}}$. For time $t\geq
t_{\D_0}$, let $\nuall^A(t)$ be the number of $r$-sets which meet
every edge of $\cH_{t}$ and contain $v$ and let $\nuall^B(t)$ be
the number of $r$-sets which meet every edge of $\cH_{t}$ and do not
contain $v$. Lemma~\ref{smdel} and $r\leq n^{5/12}$ imply that
$$ \frac{\nuall^B(t_{\D_0})}{\nuall^A(t_{\D_0})}\leq%
   n^{2/3-\Delta_0/6}. $$

Now if $t=t_{\D_0}+\sigma$ and every edge added between $t_{\D_0}$
and $t$ contains $v$ then $\nuall^A(t)\geq
\nuall^A(t_{\D_0})-\sigma$ and $\nuall^B(t)\leq
\nuall^B(t_{\D_0})$. So,
$$ \Pr(v\notin E_{t+1}\mid v\in E_j,t_{\D_0}\leq j\leq t)\leq%
   \frac{\nuall^B(t_{\D_0})}{\nuall^A(t_{\D_0})-\sigma}\leq%
   \frac{n^{2/3-\Delta_0/6}\nuall^A(t_{\D_0})}%
   {\nuall^A(t_{\D_0})-\sigma}\leq%
   2n^{2/3-\Delta_0/6}, $$
as long as $\sigma\leq \nuall^A(t_{\D_0})/2$.

Now $\nuall^A(t_{\D_0})\geq \binom{n-t_{\D_0}r}{r-t_{\D_0}}\gg
n^{\Delta_0/7}$ and so we see that if
$t_{\infty}=t_{\D_0}+n^{\Delta_0/7}$ then \whp\
\begin{equation}\label{inn}
   v\in E_t,\,t_{\D_0}\leq t\leq t_{\infty}.
\end{equation}

Let $\Omega$ denote the set of $(r-1)$-subsets of $[n]\setminus
\{v\}$ which meet $E_1,\ldots,E_{\t_4},\,\t_4=t_4-4$. (Assume a
re-numbering so that these are the edges of $\cH_{t_4}$ which do
not contain $v$). We know from \eqref{future} that
\begin{equation}\label{check}
   |\Omega|\leq (1+o(1))r^{\t_4}\binom{n}{r-\t_4}.
\end{equation}
If we condition on \eqref{inn}, the sets
$E_t\setminus\{v\},\,t_{\D_0}+1\leq t\leq t_{\infty}$ will be
chosen uniformly at random from $\Omega$ without replacement. For
a fixed $t$ and $r$-set $E\subseteq [n]\setminus \{v\}$ which meets $E_1,E_2,\ldots,E_{\t_4}$, let
$$ \pi_E=\Pr((E_{t}\setminus \{v\})\cap E\neq\emptyset) $$
and let
$$ \hat{\pi}=\max_E\{\pi_E\} . $$
Then
\begin{equation}\label{fin}
   \Pr(\exists E:\;v\notin E, E\cap E_{t}\neq \emptyset,
   t_{\D_0}+1\leq t\leq t_{\infty})
   \leq \binom{n-1}{r}\hat{\pi}^{n^{\Delta_0/7}}.
\end{equation}
(Without replacement there are fewer choices that will meet $E$).

Now fix $E$ and let $a_i=|E\cap E_i|,\,i=1,2,\ldots,t$ where
$a_1\geq a_2\geq\cdots\geq a_t$.  Also let $s=\max\{j:\;a_j\geq
r/3\}$ and note that $s\in\{0,1,2,3\}$. Then if
$\Omega_E=\{Y\in\Omega:\;Y\cap E=\emptyset\}$,
\begin{eqnarray*}
   |\Omega_E| & \geq & \brac{\prod_{j=s+1}^t(r-a_j-j)}%
   (n-2r)_{r-t}((r-t)!)^{-1}\\
   & \geq & (1-o(1))r^{t-s}n^{r-t}((r-t)!)^{-1}%
   \exp\left\{-2r^{-1}\sum a_j\right\} \\
   & \geq & (1-o(1))e^{-2}r^{t-s}n^{r-t}((r-t)!)^{-1} \\
   & \geq & \frac{|\Omega|}{10r^s}
\end{eqnarray*}
Thus,
$$ \pi_E\leq 1-\frac{1}{10r^3} . $$
It follows from \eqref{fin} that
$$ \Pr(\exists E:\;v\notin E, E\cap E_{t}\neq \emptyset,
   t_{\D_0}+1\leq t\leq t_{\infty})\leq%
   \binom{n-1}{r-1}\brac{1-\frac{1}{10r^3}}^{n^{\Delta_0/7}}=o(1) $$
and this together with Claim~\ref{schol1} finishes the proof of
Theorem \ref{th2} (except for \eqref{exactH}) since now we see
that \whp\ every edge chosen from time $t_\infty$ onwards will
contain $v$.

To prove \eqref{exactH} we re-label the edges of $\cH_{t_4}$ as
$E_1,E_2,\ldots,E_{t_4}$ so that $v\in
E_{t_4-3},E_{t_4-2},E_{t_4-1},E_{t_4}$. Recall $\t_4=t_4-4$.
 Then for $S\subseteq [\t_4]$ we let
$$ \cE_S=\left\{F\in \binom{[n]\setminus \{v\}}{r-1}:\;%
   F\cap E_j=\emptyset, j\in S\right\}. $$
Since $\cH$ is simple, if $|S|=i$ then we have
$|\cE_S|=\binom{n-ri+i(i-1)/2}{r-1}$ and \eqref{exactH} follows
directly from the inclusion-exclusion formula.

To obtain \eqref{approxH} we use \eqref{lem6a} and \eqref{t4} with
$t=t_4$ and $t_{\Delta}=\t_4$ (after observing that \whp\ from
$t_4$ on, all the edges added contain $v$). We have two asymptotic
expressions:
\begin{eqnarray*}
\frac{r^{\t_4}n^{r-\t_4-1}}{(r-\t_4-1)!}&&\\
\noalign{and}\\
\bfrac{r^2}{n}^{t_4-4}\binom{n-1}{r-1}=
\bfrac{r^2}{n}^{\t_4}\binom{n-1}{r-1}&\sim&\frac{r^{2\t_4}n^{r-1}}{n^{\t_4}(r-1)!}\qquad\qquad
\mbox{using}\ r^2=o(n) .
\end{eqnarray*}
Thus the ratio of these two expressions is asymptotically
$$\frac{(r-1)!}{r^{\t_4}(r-\t_4-1)!}=\frac{r^{\t_4}}{(r-1)_{\t_4}}\sim 1$$
because $\t_4^2/r=O(r/n^{2/3})=o(1)$. \proofend

\section{Functionals}\label{appA}
\subsection{Graph functionals} We introduce a certain type of
graph functional. Here $G$ is a graph with $s$ vertices and $f$
edges. $\cM$ denotes set of graphs which are matchings.

\begin{proposition}\label{Gfun}
Let $t$ be a positive integer and $x$ and $y$ be nonnegative
quantities such that $x=o(1)$, $y=\omega(1)$ and $t^2xy^2=O(1)$.
Then,
\begin{eqnarray}
   \sum_{\emptyset\neq G\not\in\cM}x^fy^s
   &=& O(t^4x^2y^3)\label{termerr}
   \\
   \sum_{G\in \cM,f\geq 2}x^fy^s
  &=&O(t^4x^2y^4)\label{estf2} \\
   \sum_{G\in \cM,f=1}x^fy^s
   & = & \bin{t}{2}xy^2 .\label{termf2}
\end{eqnarray}
\end{proposition}
\proofstart
Let $n_t(f,s)$ count the number of
subgraphs of $K_{[t]}$ that have exactly $f$ edges and exactly
$t-s$ isolated vertices. $n_t(f,s)\leq\bin{\bin{t}{2}}{f}$.  Let $f_0=\lceil
(s+1)/2\rceil$.
\begin{eqnarray*}
 \sum_{\emptyset\neq G\not\in\cM}x^fy^s
   & = & \sum_{s=3}^t\sum_{f=f_0}^{\bin{s}{2}}
         n_t(f,s)x^fy^s \\
   & \leq & \sum_{s=3}^t\sum_{f=f_0}^{\bin{s}{2}}
            \frac{\bin{t}{2}^{f}}{f!}x^fy^s \\
   & \leq & \sum_{s=3}^ty^s
            \frac{(t^2x/2)^{f_0}}{f_0!}
            \sum_{f=0}^{\bin{s}{2}-f_0}
            \frac{\bin{t}{2}^f}{f!}x^f \\
   & \leq & \sum_{s=3}^ty^s
            \frac{(t^2x/2)^{f_0}}{f_0!}
            \sum_{f=0}^{\binom{s}{2}-f_0}
            \frac{1}{f!}\left(\frac{t^2x}{2}\right)^f \\
   & \leq & \exp\left\{\frac{t^2x}{2}\right\}\brac{
            \sum_{t=0}^{\lfloor t/2\rfloor-2}
            \frac{(t^2x/2)^{t+2}y^{2t+3}}{t!}+
        \sum_{t=0}^{\lfloor t/2\rfloor-2}
        \frac{(t^2x/2)^{t+3}y^{2t+4}}{t!}} \\
   & = & (1+o(1))(t^4x^2y^3/4+t^6x^3y^4/8)\sum_{t=0}^{\lfloor t/2\rfloor-2}
\frac{(t^2xy^2/2)^t}{t!}   \\
&\leq&t^4x^2y^3\exp\left\{\frac{t^2xy^2}{2}\right\}\\
&=&O(t^4x^2y^3).
\end{eqnarray*}

In the case where $F$ is a matching,
$n_t(f,s)=\frac{(t)_{2f}}{2^f}\frac{1}{f!}$ because $s=2f$.
Thus
\begin{eqnarray*}
 \sum_{G\in \cM,f\geq 2}x^fy^s   & = & \sum_{f=2}^{\lfloor t/2\rfloor}
                        \frac{(t)_{2f}}{2^f}\frac{1}{f!}
                        x^fy^{2f} \\
                  & \leq & \frac{1}{2}
                           \left(\frac{t^2xy^2}{2}\right)^2
                           \exp\left\{\frac{t^2xy^2}{2}\right\}
\end{eqnarray*}
and \eqref{estf2} follows. Equation \eqref{termf2} is clear.
\proofend

\subsection{Grid functionals}\label{gridfun}

Now we introduce a functional on a grid.

\begin{proposition}\label{16}
Suppose that
$\t$ is a positive integer and that $x,y$ are positive reals such that
$$x=o(1),\,y=\omega(1),\t\Delta xy^2=o(1).$$
Then
$$\sum_{\substack{\emptyset\neq H\subseteq [\t]\times [\Delta] \\
 u=\left|H\mid_{[\t]}\right| \\
\ell=\left|H\mid_{[\Delta]}\right| }}
      x^h y^{u+\ell}=(1+o(1))\t\Delta xy^2.$$
\end{proposition}
\proofstart
Let $n(h,u,\ell)$ denote the number of sets $H\subseteq
u\times\ell$ with $|H|=h$, $\left|\left. H\right|_{[u]}\right|=u$
$\left|\left. H\right|_{[\ell]}\right|=\ell$
$$ \sum_{\substack{
         H\subseteq [\t]\times [\Delta] \\
         u=\left|H\mid_{[\t]}\right| \\
         \ell=\left|H\mid_{[\Delta]}\right|
         }}x^hy^{u+\ell}=
   \sum_{h,u,\ell}\binom{\t}{u}\binom{\Delta}{\ell}
   n(h,u,\ell)x^hy^{u+\ell} $$

We use the bounds
$$ n(h,u,\ell)\leq\left\{\begin{array}{ll}
                         \ell^u
                         \binom{u\ell}{h-\ell}, & \mbox{ if
                         $u\geq\ell$;} \\
                         u^\ell
                         \binom{\ell u}{h-u}, & \mbox{ if
                         $\ell\geq u$.}\end{array}\right. $$
Indeed, if $u\geq\ell$ then
\begin{eqnarray*}
n(h,u,\ell)&\leq&\sum_{t=1}^\ell\binom{u}{t}\sum_{\substack{d_1+\cdots+d_t=\ell\\
d_1\geq 1,\ldots,d_t\geq 1}}\frac{\ell!}{d_1!\cdots d_t!}\ell^{u-t}\binom{u\ell-\ell-u+t}{h-\ell-u+t}\\
&\leq&\sum_{t=1}^\ell\binom{u}{t}\ell^{u-t}t\ell\binom{u\ell-\ell-u+t}{h-\ell-u+t}\\
&\leq&\sum_{t=1}^\ell\binom{u}{u-t}\ell^u\binom{u\ell-u}{h-\ell-u+t}\\
&\leq&\ell^u\binom{u\ell}{h-\ell}\qquad\mbox{by the Vandermonde identity.}
\end{eqnarray*}
Similarly, if $\ell\geq u$ then $n(h,u,\ell)\geq
u^\ell\binom{u\ell}{h-u}$.

Therefore, we can bound the summation by four other summations.

\begin{eqnarray}
   \lefteqn{\sum_{\substack{
                  \emptyset\neq H\subseteq [\t]\times [\Delta] \\
                  u=\left|H\mid_{[\t]}\right| \\
                  \ell=\left|H\mid_{[\Delta]}\right|
                  }}x^hy^{u+\ell}} \nonumber \\
   & \leq & \sum_{\ell\geq 1}\binom{\t}{\ell}\binom{\Delta}{\ell}
            \ell^\ell\left(xy^2\right)^{\ell} \qquad Case:\ h=u=\ell\label{heql} \\
   & & +\sum_{\ell\geq 2}\sum_{h>\ell}\binom{\t}{\ell}
        \binom{\Delta}{\ell}\ell^\ell\binom{\ell^2}{h-\ell}
        x^hy^{2\ell} \qquad Case:\ h>u=\ell\label{hgtl} \\
   & & +\sum_{\ell\geq 1}\sum_{u>\ell}\sum_{h\geq u}
        \binom{\t}{u}\binom{\Delta}{\ell}
        \ell^u\binom{u\ell}{h-\ell}x^hy^{u+\ell}\qquad Case:\ h\geq u>\ell
        \label{ugtl} \\
   & & +\sum_{u\geq 1}\sum_{\ell>u}\sum_{h\geq\ell}
        \binom{\t}{u}\binom{\Delta}{\ell}
        u^\ell\binom{\ell u}{h-u}x^hy^{u+\ell}\qquad Case:\ h\geq\ell>u
        \label{lgtu}
\end{eqnarray}
Consecutive terms in (\ref{heql}) increase by a factor of size $O(\t\Delta xy^2)=o(1)$ and so
the sum is dominated by the first term i.e.
$$RHS(\ref{heql})=(1+o(1))\t\Delta xy^2.$$
Let $m=\min\{\D,\t\}$.
\begin{eqnarray*}
RHS(\ref{hgtl})&\leq&\sum_{\ell=1}^m\binom{\t}{\ell}
            \binom{\Delta}{\ell}\ell^\ell y^{2\ell}
            \sum_{h'\geq 0}\binom{\ell^2}{h'+1}
            x^{h'+\ell+1} \\
   & \leq &(1+o(1))\sum_{\ell=1}^m\binom{\t}{\ell}
            \binom{\Delta}{\ell}
            \left(\ell xy^2\right)^{\ell}\ell^2x
            e^{\ell^2x} \\
   & = & (1+o(1))\t\Delta x^2y^2.
\end{eqnarray*}
since successive terms increase by  a factor of size $O(\t\Delta xy^2)$.

Now we bound summation (\ref{ugtl}).

\begin{eqnarray*}
RHS(\ref{ugtl})
   & \leq & \sum_{\ell\geq 1}\sum_{u>\ell}\sum_{h\geq u}
        \binom{\Delta}{\ell}\binom{\t}{u}
        \ell^u\binom{u\ell}{h-u}x^hy^{u+\ell}\\
& \leq & \sum_{\ell\geq 1}\sum_{u>\ell}
         \binom{\Delta}{\ell}\binom{\t}{u}({\ell}xy^2)^u
         \sum_{h'\geq 0}\binom{u\ell}{h'}x^{h'} \\
   & \leq & (1+o(1))\sum_{\ell\geq 1}\sum_{u>\ell}
          \binom{\Delta}{\ell}  \binom{\t}{u}({\ell}xy^2)^u\\
   & = & (1+o(1))\sum_{\ell\geq 1}\binom{\Delta}{\ell}
\binom{\t}{\ell+1}(\ell xy^2)^{\ell+1}\\
   & \leq &(1+o(1))\t^2\Delta x^2y^4
\end{eqnarray*}

By symmetry, (\ref{lgtu}) is  also bounded by $(1+o(1))\t\Delta^2 x^2y^4$.

The sum of \eqref{heql}, (\ref{hgtl}), (\ref{ugtl}) and (\ref{lgtu}) is therefore $(1+o(1))\t\Delta xy^2$.
\proofend

\section{Remarks}

The thresholds, $r=\Theta(n^{1/3})$ and $r=\Theta(n^{5/12})$ each
present their own unique difficulties.  When $r=cn^{5/12}$ there
is a probability, $p_1=p_1(c)>0$ that the first intersection of
size greater than $1$ will occur before a degree three vertex.  In
our notation, $\tbigint<t_3$.  Our analysis will not work if this
occurs.

The threshold $r=\Theta(n^{1/3})$ presents a different problem.
When the first degree three vertex emerges, it may not be unique.
For example, if $r=cn^{1/3}$, then there is a probability
$p_2=p_2(c)>0$ that $\{e_1,e_2,e_3,e_4\}$ form a simple hypergraph
but $e_5$ contains both $e_1\cap e_2$ and $e_3\cap e_4$.  What
makes this case even more difficult is that $t_4$ need not be
$t_3+1$.  Recall that we proved that when
$\omega(n^{1/3})=r=o(n^{2/5})$ then $t_4=t_3+1$.  At
$r=\Theta(n^{1/3})$ this is not the case.  For example, if
$r=cn^{1/3}$, there exists a probability $p_3=p_3(c)>0$ such that
all of the following occurs: The edges $\{e_1,\ldots,e_6\}$ form a
simple hypergraph.  The edge $e_7$ contains intersection points
$e_1\cap e_2$ and $e_3\cap e_4$ and no others.  The edge $e_8$
contains intersection points $e_1\cap e_2\cap e_7$ and $e_5\cap
e_6$ but no others. Then, $e_9$ contains intersection points
$e_3\cap e_4\cap e_7$ and $e_5\cap e_6\cap e_8$ but no others. So,
$t_4\geq t_3+2$ for this example.

In fact, there are numerous outcomes that can occur with nonzero
probability after only $O(1)$ edges when $r=\Theta(n^{1/3})$.  The
simplicity of Theorem~\ref{th2} is, therefore, all the more
remarkable.

A subset of the authors of this paper intend to work further to
describe the hypergraph that results when $r$ is a threshold value
as well as proceed to the case where $r=\omega(n^{5/12})$.  We
believe that precise structural results such as
Theorems~\ref{main} and~\ref{th2} are impossible, but the size of
the hypergraph may be able to be described.

\section{Thanks}
We would like to thank some anonymous referees for careful reading
and useful suggestions to improve the manuscript.

\end{document}